\documentclass[a4paper, reqno]{amsart}
\usepackage[T1]{fontenc}
\usepackage{a4wide}

\usepackage{amsmath}
\usepackage{amssymb}
\usepackage{enumerate}	
\usepackage{nicefrac}		
\usepackage{array}		
\usepackage{calc}			
\usepackage{flafter}		
\usepackage{esdiff}
\usepackage{pgffor}   
\usepackage{pgf} 
\usepackage{mathabx}
\usepackage{multicol}
\usepackage{soul} 
\definecolor{newblue}{RGB}{94,75,250}
\definecolor{newblue2}{cmyk}{1,0.6,0,0.06}
\definecolor{grey}{gray}{0.5}

\usepackage{multirow}

\usepackage{tikz}
\usetikzlibrary{calc}

\usepackage
[
pdfauthor   = {Lisa Maria Kreusser},
pdftitle    = {Title},
pdfsubject  = {},
pdfkeywords = {},
bookmarks=true,
bookmarksopen=true,
pagebackref=false,
hyperindex=true,
colorlinks=true,
linkcolor=blue,
citecolor=blue,
filecolor=blue,
urlcolor=blue,
]
{hyperref}

\usepackage{mathrsfs}	
\usepackage{graphicx}
\usepackage[caption=false,subrefformat=parens,labelformat=parens]{subfig}

\usepackage[nameinlink]{cleveref}
\usepackage{verbatim}

\usepackage{bbm}

\newcommand\N{\mathbb{N}}
\newcommand\R{\mathbb{R}}

\newcommand\bl{\left(}
\newcommand\br{\right)}
\newcommand\Vbd{\partial V}

\newcommand*\di{\mathop{}\!\mathrm{d}}

\renewcommand\epsilon{\varepsilon}
\renewcommand\theta{\vartheta}

\newtheoremstyle{mytheoremstyle} 
{6pt}                   
{6pt}                   
{\itshape}                  
{}					
{\bf}                 
{}                         
{.5em}                       
{}  
\theoremstyle{mytheoremstyle}
\newtheorem{satz}{Satz}[section]

\theoremstyle{remark}

\newtheorem{remark}[satz]{Remark}

\graphicspath{{./../}}

\numberwithin{equation}{section}
\usepackage{cleveref}

\title{Models for information propagation on graphs}
\author{Oliver R.\ A.\ Dunbar, Charles M.\ Elliott and Lisa Maria Kreusser}

\begin{document}

\maketitle

\begin{abstract}
We propose and unify classes of different models for information propagation over graphs. In a first class, propagation is modelled as a wave which emanates from a set of \emph{known} nodes at an initial time, to all other \emph{unknown} nodes at later times with  an ordering determined by the arrival time of the information wave front. A second class of models is based on the notion of a travel time along paths between nodes. The time of information propagation from an initial  \emph{known} set of nodes to a node is defined as the minimum of a generalised travel time over subsets of all admissible paths. A final class is given by imposing a local equation of an eikonal form at each \emph{unknown} node, with  boundary conditions at the \emph{known} nodes. The solution value of the local equation at a node is coupled to those of neighbouring nodes with lower values. We provide  precise  formulations of  the model classes and prove equivalences between them. Finally we apply the  front propagation models on graphs to semi-supervised learning via label propagation and information propagation  on trust networks.
\end{abstract}

\section{Introduction}
Information propagation (also known as diffusion, cascade, or spread) is of great importance in complex networks where, given  information at a small number of nodes of the network, the aim is to understand the  propagation to all the nodes. Social media networks provide typical examples including  the breaking of a news story and  the spread of product advertisements,  internet memes and  misinformation to different users. The ability to predict propagation plays a key role in tasks such as informing how to seed information for obtaining maximal coverage and influence \cite{KemKleTar15,JalPer17}, or for identifying likely sources of information provided that the times are given when the information was received \cite{Bol_etal18}. Models may be used for control and management of the propagation.  

Our starting point is to model the  underlying network as a  given graph. The aim of this work is to formulate models for inspired by the
propagation  of waves passing through continuous media.  Elements of the approach are  that information has either arrived at a graph vertex or not, that information is transmitted to a node only from neighbouring nodes at which information has arrived already and that there is an arrival time for each  node. These models for information propagation can then be used in applications ranging from social media networks to semi-supervised learning.

\subsection{Continuum front propagation}
In the  continuum setting,  there are three common viewpoints for modelling waves: front propagation, first arrival times and local equations. To introduce these viewpoints, we consider an open bounded domain $\Omega \subset \R^d$ for $d\geq 1$ with a Lipschitz boundary $\Gamma$, a given  point $x_0 \in \Omega$ and a continuous, positive function $s  \colon \bar\Omega\to \R$ which can be regarded as the impedance of the medium $\bar \Omega$.

A first approach proposes a propagating front separating the region for which the wave has arrived from the remainder. The fronts initiate at $x_0$, and are characterised by being level surfaces of the arrival time from $x_0$. The impedance $s(x)$ is specific for the underlying medium and controls the additional time required for the front to travel through the medium at $x$. We also refer to this approach as front propagation. 

A second classical approach consists of formulating a model based on finding the smallest travel time over a set of possible paths and hence results in an optimisation problem. The aim of this model is to determine the shortest travel time along any path from $x_0$ to every $x \in\bar\Omega, x\neq x_0$, in the medium $\bar \Omega$ for a given impedance $s$.  
This task can be expressed as the minimisation problem 
\begin{equation}\label{eq:eikonalmin}
	u(x)=	\inf_{\substack{\xi \in W^{1,\infty}([0,1],\bar \Omega),\\\xi(0)=x_0,\enspace\xi(1)=x}}\left\{ \int_0^1 s(\xi(r)) \|\xi'(r)\|_2\di r \right\},
\end{equation}
cf.\ \cite{DecEllSty11}, where $\|\cdot \|_2$ denotes the 2-norm in $\R^d$ and $\xi(\cdot)$ is a parameterised path in the Sobolev space $W^{1,\infty}$. Note that $\xi \in W^{1,\infty}([0,1])$ is locally Lipschitz continuous and hence the integral in \eqref{eq:eikonalmin} is well-defined. Since large values of $s$ slow down the movement and increase  the travel time within the medium, we sometimes refer to $s$ as the slowness function, while its inverse  $\tfrac{1}{s}$ can be regarded as a velocity. We also refer to this approach as first arrival times.

A third approach arises when regarding an optimal value $u$ of \eqref{eq:eikonalmin}  as  a solution to the eikonal equation, an isotropic static Hamilton-Jacobi partial differential equation. The eikonal equation is given by
\begin{equation}\label{eq:eikonal}
	\|\nabla u\|_2=s \quad \text{in} \enspace\Omega \backslash \{x_0\}
\end{equation}
with boundary conditions 
\begin{align}\label{eq:eikonalbc}
\begin{split}
	u(x_0)&=0,\\
	\nabla u(x)\cdot \nu(x)&\geq 0 \quad \text{for}\enspace x\in \Gamma,
\end{split}
\end{align}
where $\nu$ is the unit outer normal to $\Gamma$. We also refer to this approach as  a local equation. Also it is possible
to pose and solve eikonal equations on connected (sub)Riemannian manifolds, see e.g.\
\cite{Gromov2007}.

These three approaches of wave propagation in continuum settings have been exploited to advance different fields of research. The optimisation over paths (also referred to as first arrival times)  arises in  modelling of optimal logistics such as accessibility, evacuation planning, robot navigation and ray models. The study of the graph  eikonal equation (i.e.\ a local equation) is of importance for proving  theoretical results on existence and uniqueness of solutions with certain monotonicity properties. Efficient numerical methods such as fast marching algorithms  take advantage of the front propagation approach when solving the continuum eikonal equation  \cite{Set96,KimSet98,Set99}. This demonstrates that diverse perspectives on modelling waves are crucial in the continuum setting for getting more insights into modelling, analysis and numerical methods of the underlying continuum problem.

In contrast to the continuum setting, only a scattered picture is currently available for graphs, including shortest paths, Dijkstra’s algorithm and graph-eikonal models.
Motivated by the continuum setting, the aim of this work is to propose and  unify  corresponding perspectives in the graph setting.  We formulate and relate  several  classes of models based on front propagation,  first arrival time over sets of admissible paths and  a local equation considering arrival times at a given node and its neighbours. As part of this, we introduce appropriate graph-based generalisations of the continuum counterparts for the three classes of models. In  the context of the Dijkstra algorithm, for instance, the Dijkstra algorithm can be regarded as a front propagation model. For the local equation, we replace \eqref{eq:eikonal}--\eqref{eq:eikonalbc} in the continuum setting by a graph-based version of the  local equation 
\[
	\|\nabla u\|_p=s \quad \text{in} \enspace\Omega \backslash \{x_0\}
\]
for $p=\infty$ with  boundary conditions \eqref{eq:eikonalbc}, which leads to an $\ell^\infty$ graph-eikonal equation. We also propose a first arrival time model, based on the travel time over paths, and prove its equivalence to Dijkstra's algorithm. Motivated by the special case $p=\infty$ for the local equation, we derive front propagation, first arrival time and local equations for other cases of $p$. The main contribution of this paper is to model wave propagation in the graph-based setting using three perspectives (front propagation, first arrival times and local equations). We prove the equivalence of the models for  special cases of  $p$. It is important to note that  in the models we do not embed the vertices in any ambient Euclidean space. 

\subsection{Applications}
 It is natural to introduce the concept of information propagation to data classification and semi-supervised learning.  Motivated by this, we apply front propagation on graphs to  classical examples  in semi-supervised learning such as the the {\it Two moons problem} and  { \it Text classification datasets}. Here the information consists of a given finite set of labels  and the aim is to    label all vertices in a graph based on the  knowledge of the labels on given small number of  nodes. Labels are attached by ordering the magnitudes of the arrival times of the information.  In addition, we  apply information propagation to {\it Trust networks}. These are  social networks whose  users rate each other by trustworthiness.  Examples include collaborative networks such as a community of software engineers, or partners of a transaction within cryptocurrency exchanges. Applied to the software community dataset \texttt{soc-advogato} \cite{AhmRos}, we show that information propagation can use local trust information to create rankings of any collaborator on the network. Our model-rankings are resistant to Sybil attack \cite{Dou02,LuSheZha14,AlQ_etal17}, where users artificially inflate their reputation, by creating a group of fake users to giving them positive ratings.

\subsection{PDE approaches}
 Many computational methods for semi-supervised and unsupervised classification \cite{BluCha01, BelMatNiy04, Zhu05} are based on variational models and PDEs \cite{SleTri16}.  Examples include algorithms  based on phase fields \cite{ber12} and the MBO scheme \cite{mer13}, as well as $p$-Laplacian equations \cite{Elm15,Kre20}.  
 In a series of papers, Elmoataz et al.\ \cite{DesElmLez13,DesElmLezTa10,BouElmLez08, ElmLezTa09} postulate discrete eikonal equations and investigate label propagation on graphs with applications in imaging and machine learning. Current analytical results include  an investigation of viscosity solutions for Hamilton-Jacobi equations on networks \cite{CamMar13}, the well-posedness of nonlinear PDEs such as the Eikonal equation on finite graphs \cite{Oberman2015} and an approximation scheme for an eikonal equation on a network  \cite{CamFesSch13}, producing an approximation of shortest paths to the boundary. In addition, limits and consistency of non-local and graph approximations to  the time-dependent (local) eikonal equation have been studied in \cite{Fad21}. The robustness of the solution to the eikonal equation for $p=1$ and its convergence to the shortest path distance as $p\to \infty$ is shown in \cite{CalEtt22}.

\subsection{Contributions}
Our contributions are as follows:
\begin{itemize}
    \item Derivation of  general model formulations for three perspectives (front propagation models, first arrival time models and local equations) in the graph-based setting that include established models (Dijkstra's algorithm, shortest paths and $\ell^p$ graph-eikonal equations for $p\geq 1$) as special instances.
    \item Unification of the three perspectives in the graph-based setting by proving equivalence of the models (front propagation, first arrival times and discrete generalised eikonal models) depending on $p$.
    \item Application of front propagation on a weighted social network to calculate metrics of trust securely.
    \item Application of front propagation on graphs to classical problems in semi-supervised learning for point cloud data sets (two moons problem, text classification datasets Cora and CiteSeer).
\end{itemize}

\subsection{Outline}
We introduce several  models for travel times on a graph in Section~\ref{sec:models}. Equivalences between certain instances of the models are established in  Section \ref{sec:updateformula}.  
In Section \ref{sec:trustnetworks}, we apply information propagation to trust networks. The use of front propagation on graphs to semi-supervised learning via label propagation is illustrated in Section \ref{SemiSuper}. Finally we make some concluding remarks in Section \ref{Conclusion}.

\subsection{Notation}\label{sec:notation}
Following  the terminology and setting in \cite{DunSleStuTho18_pre,BouElmLez08, SleTri16}, we consider a finite, 
connected weighted graph $G=(V,E,w)$   with vertices $V=\{1,\ldots,n\},$ 
edges $E\subset V^2$ and  nonnegative edge weights $w$. We assume that the graph is simple, i.e.\ there exists at most one edge between any two vertices. 
We suppose that there is a decomposition of $V:=\partial V\cup \mathring V$ into two disjoint non-empty sets $\partial V$ and $\mathring V$.
The edge between node $i$ and node $j$ is denoted by $(i,j)$.  For ease of notation, we regard the weights $w$ as a weight matrix $w\in \R^{n\times n}$ with entries $w_{ij}$, where we assume that there exists an edge $(i,j)\in E$ if and only if $w_{ij}>0$, while $w_{ij}=0$ if $(i,j)\notin E$.
 Since $G$ is not necessarily undirected, $w_{ij}\neq w_{ji}$ in general.
This framework also includes unweighted graphs corresponding to the cases in which $w_{ij}=1$ for all $(i,j)\in E$. 
Given a graph $G$, we denote by  $N(i)\subset V$  the set of neighbours of node $i\in V$ 
. We define $j\in N(i)$ if there exists an edge $(j,i)\in E$, and in general this does not imply existence of $(i,j) \in E$. The direction of this relationship is chosen for convenient notation in the following.
We introduce the notion of a path from node $x\in V$ to $y\in V$ and write  $p_{x,y}=(x=i_1,\ldots,y=i_{n(p_{x,y})})$ for a path with $n(p_{x,y})$ nodes and $n(p_{x,y})-1$ edges $(i_{m-1},i_m)\in E$ for $m=2,\ldots,n(p_{x,y})$ such that all nodes $i_m$ for $m\in\{1,\ldots,n(p_{x,y})\}$ are  distinct, i.e.\ a path must not self-intersect.
Due to the assumption that the graph $G$ is connected, for every $x,y \in V$ there exists a path $p_{x,y}$  connecting $x$ and $y$, i.e.\ there exists  $n(p_{x,y})>1$ such that  $p_{x,y}=(x=i_1,\ldots,y=i_{n(p_{x,y})})$ is a path with edges $(i_{m-1},i_m)\in E$ for $m=2,\ldots,n(p_{x,y})$.  
For a graph with $|V|=n$ nodes, we denote by $\mathcal H^n$ the function space of all functions defined on $V$, i.e.\ all $v\in \mathcal H^n$ are of the form $v\colon V\rightarrow \mathbb R$. For $v\in \mathcal H^n$, we  write $v_x=v(x)$ for $x\in V$. 
We also assume that  there is a given {\it slowness} function   $s\in \mathcal H^n$ with $s\ge 0$.

\section{Description of  models}\label{sec:models}
In this section, we propose several  models for the propagation of information on graphs. The common elements of the models are:
\begin{itemize}
\item We suppose that either all information has  arrived at a vertex or none.
\item
We introduce the variable  $u\in \mathcal H^n$ with  $u_i$ for  $i\in V$ to denote the {\it arrival time} of information at vertex $i$.
\item We assume that $u$ is prescribed on $\partial V$ and we set $u=0$ on $\partial V$, though in general the models can accommodate a wider class of boundary conditions.
\item
We suppose that information propagation is local. That is, information arrives at a vertex only by propagation from a neighbouring vertex for which information has arrived.  Thus  there is a unique  travel time $u_i$ at each  node $i$ that can only depend on travel times at nodes $j\in N(i)$ with $u_j < u_i$. 

\item
The edge weights reflect the distance or resistance to propagation along an edge.
\item
The function $s\in \mathcal H^n$ is a measure of slowness or resistance  associated with each vertex.

\end{itemize}
 The aim of a model is to associate a travel time $u_i$ with each vertex of the graph. Since the graph is finite,  $u=\{u_i,i=1,2,\ldots,n\}$  attains an unknown number of $J+1\in \N$ distinct   values 
consisting of prescribed initial data $U_0\in \mathbb R$ and unknown values $U_1,\ldots, U_J\in \R$  ordered so that $U_0<\ldots < U_J$.
We set $V_0:=\Vbd$ as the set of initially labelled vertices and prescribe the initial data $U_0$, i.e.\ $u_i=U_0$ for all $i\in V_0$. In the following we set $U_0=0$. 
We consider three classes of models. 
The first class of models is based on the propagation of discrete fronts from an initial front $\partial V$ (Model 1).  The second  class of models considers  first arrival  times of sets of paths that link vertices in the initial set $\partial V$ to vertices in $\mathring V=V\backslash \Vbd$ (Model 2). For the third class of models, we postulate a generalised discrete $\ell^p$ eikonal equation model (Model 3) depending on parameter $p$. We mainly focus on $p\in\{1,2,\infty\}$ below. Note that some of the model instances may look rather complicated. However, the main motivation is to unify graph-based models from three perspectives (front propagation, first arrival times, local equations) by proving their equivalence.

\subsection{Front propagation models}\label{sec:frontpropagationmodels}

In this approach, we view information propagation   as an evolving front, i.e.\ a boundary that separates the region for which the wave has arrived from the remainder. We decompose the set $\mathring V$ of initially unlabelled vertices into $J$ disjoint sets $V_1,\ldots,V_J$ such that for $j\in \{1,\ldots, J\}$ all vertices $i\in V_j$
satisfy $u_i=U_j$.
We define  \emph{known} sets  $K_0,\ldots,K_J$ and  \emph{candidate} sets $C_0,\ldots, C_J$ as follows:
\[
K_l = \bigcup_{j \in\{0,\dots,l\}}\ V_j, \qquad  C_l = \bigcup_{j \in K_l}\  N(j) \ \setminus \ K_l.
\]
Under the assumption  that  $U_j$ and $V_j$ for $j=0,\ldots,k-1$ are known, implying that the value of $u_i$  for all $i \in K_{k-1}$ is known, our task is  to determine $U_k$ and $V_k$. The front $F_{k-1}$ consists of all vertices in $K_{k-1}$ with neighbours in $C_{k-1}$ and with $F_0=V_0$. 
We determine candidate values  $\tilde{u}_i$ for each $i\in C_{k-1}$ using a model (specified below) and we define $U_k$ by choosing the smallest candidate value  in  the candidate  set $C_{k-1}$:
\begin{equation}\label{eq:propagationmodel}
U_k := \min_{i\in C_{k-1}} \tilde{u}_i.
\end{equation}
We then define $V_k\subset C_{k-1}$ to be the set where the minimum  is attained and we set $u_i=U_k$ for all $i\in V_k$. 
The above procedure depends on the definition of candidate values $\tilde u_i$ for $i\in C_{k-1}$. We define relationships for $\tilde u_i$  that depend upon the set $N(i)\cap K_{k-1}$. Using \eqref{eq:propagationmodel}, the values $U_1,\ldots,U_L$ of the solution $u$ can then be determined. By construction, the solution $u$ is unique for the function $i\mapsto \tilde u_i$.

\subsubsection{Model 1(i)}
Given the known arrival time $u_j$ for  $j\in K_{k-1}$, and let $j\in N(i)$ so edge $(j,i)$ exists, then a candidate for the arrival time at $i$,  is given by $u_j + \tfrac{s_{i}}{w_{j,i}}$. Choosing the smallest value of all these possible candidate values results in the candidate 
\begin{equation}\label{eq:frontmodel_i}
\tilde{u}_i =\min_{j\in N(i)\cap K_{k-1}} \left\{u_j + \frac{s_{i}}{w_{j,i}}\right\} 
\end{equation}
for $i\in C_{k-1}$. 
Here,  $u_j + \tfrac{s_{i}}{w_{j,i}}$ is the sum of the first arrival time $u_j$   at node $j$ and $\tfrac{s_{i}}{w_{j,i}}$ which is the travel time from $j$ to $i$ along edge $(j,i)$. The travel time along $(j,i)$ only depends on the slowness $s_i$  at the endpoint of $(j,i)$ and the edge weight $w_{j,i}$. The term $\tfrac{s_{i}}{w_{j,i}}$ is inspired from the continuum setting \eqref{eq:eikonalmin} which suggests that the travel time along an edge $(i,j)$ is antiproportional to the velocity $\tfrac{1}{s_i}$ and hence proportional to $s_i$. \eqref{eq:eikonalmin} also suggests that the travel time is proportional to the length of an edge and thus proportional to $\tfrac{1}{w_{i,j}}$ if we regard $w_{i,j}$ as a characterisation of the connectivity of vertices $i$ and $j$. 

As the minimum in \eqref{eq:frontmodel_i} can be associated with the $\ell^\infty$-norm, we will also see later that this model is equivalent to the $\ell^\infty$ graph-eikonal equation.

\subsubsection{Model 1(ii)}\label{sec:model2ii}
While only the smallest neighbouring value has been considered in \eqref{eq:frontmodel_i} which can be associated with the $\ell^\infty$-norm, we consider a more averaging approach in the following instance of a front propagation model motivated by weighing neighbouring known values in an $\ell^2$-sense.
We define
$z_i^2:= \sum_{j\in N(i)\cap K_{k-1}} w_{j,i}^2$ for $i\in C_{k-1}$, i.e.\ 
$z_i^2=\|(w_{j,i})_{j\in N(i)\cap K_{k-1}}\|_2^2$. 
For $i\in C_{k-1}$, we set
\begin{equation}\label{eq:frontmodel_ii}
\tilde u_i=\mu_i+ \sqrt{\frac{s_i^2}{ z_i^2}  -\sigma^2_i}.
\end{equation}
Here, 
\[\mu_i=\frac{1}{ z_i^2} \sum_{j\in{N}(i)\cap K_{k-1}} w_{j,i}^2u_j\]
can be regarded as the weighted mean travel time between any node $j\in {N}(i)\cap K_{k-1}$ and node $i$ as $\frac{1}{ z_i^2} \sum_{j\in{N}(i)\cap K_{k-1}} w_{j,i}^2=1$. The weighted mean travel time to $i$ balances the travel time to each  known node $j$ with the squared weights between $i$ and $j$. Further, we set
\[\sigma^2_i=\sum_{j\in {N}(i)\cap K_{k-1}}\left( \frac{w_{ i,j}^2}{z_i^2} u_j^2\right)-\mu_i^2\]
as the variance of the weighted mean travel time. 

As an interpretation of \eqref{eq:frontmodel_ii}, we can regard the wavefront of information travelling simultaneously from all known nodes $j\in K_{k-1}$ to candidate node $i$ where the averaged wavefront (in the $\ell^2$-sense) depends on the weighted mean travel time $\mu_i$ and its variance  $\sigma_i^2$.  With this model, one can interpret the neighbours' values as forming an estimate of a candidate value $\tilde{u}_i$ from below, with a weighted mean square error $(\tilde{u}_i - \mu_i)^2 + \sigma^2_i = \frac{s^2_i}{z^2_i}$. We will also see later that this model is equivalent to the $\ell^2$ graph-eikonal equation.

 \subsubsection{Model 1(iii)}
Similarly to \eqref{eq:frontmodel_ii}, we  consider an averaging approach in the following instance of a front propagation model, but here we weigh neighbouring known values in an $\ell^1$-sense.
For $i\in C_{k-1}$, we define $M_{i,k}=| N(i)\cap K_{k-1} |$ and  $y_i :=\sum_{j\in N(i)\cap K_{k-1} }w_{j,i}$, i.e.\ $y_i=\|(w_{j,i})_{j\in N(i)\cap K_{k-1}}\|_1$. We set
\begin{equation}\label{eq:frontmodel_iii}
\tilde{u}_i  =\frac{1 }{y_i} \sum_{j\in N(i)\cap K_{k-1}} (w_{j,i}u_j)+\frac{s_i}{y_i}  =\frac{1 }{y_i} \sum_{j\in N(i)\cap K_{k-1}} w_{j,i} \left(u_j+\frac{s_i}{M_{i,k}w_{j,i}}\right)
\end{equation}
for $i\in C_{k-1}$. The first term in \eqref{eq:frontmodel_iii} can be regarded as a weighted mean travel time to $i$, obtained by balancing the travel time from each known node $j$ with the weight $w_{j,i}$ between $j$ and $i$, while the second term $\tfrac{s_i}{y_i}$ can be interpreted as bias. Like for the other instances, we can interpret \eqref{eq:frontmodel_iii} as the wavefront of information travelling simultaneously from all known nodes $j\in K_{k-1}$ to candidate node $i$ where the averaged wavefront  (in the $\ell^1$-sense) depends on the weighted mean travel time  and its bias. We will also see later that this model is equivalent to the $\ell^1$ graph-eikonal equation.

\subsection{First arrival times}\label{sec:pathset}

\begin{figure}[htbp]
    \centering
    \includegraphics[width=0.8\textwidth]{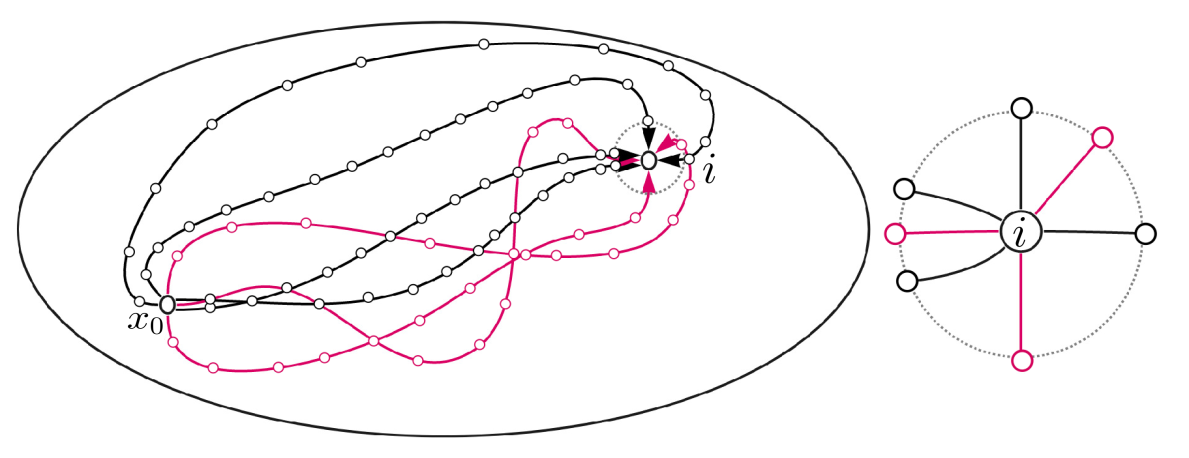}
    \caption{An illustration of a path set and its truncation. On the left we represent the set of all paths $\mathbb{P}_{x_0,i}$ between two nodes $x_0$ and $i$ with black arrows from $x_0$ to $i$. We represent a path set $P_{x_0,i} \subset \mathbb{P}_{x_0,i}$ in pink. In particular the path set $P_{x_0,i}$ contains three paths. On the right of the figure, we zoom into the neighbourhood $N(i)$, represented as nodes on dotted circle; the pink nodes on the dotted circle represent the penultimate truncation $K(P_{x_0,i})\subset N(i)$ of the path set. The pink edges therefore can be written as $(j,i)$ such that $j \in K(P_{x_0,i})$.  }
    \label{fig:pathsets}
\end{figure}

In this approach, we optimise travel times over path sets as a generalisation of travel times over paths. For this, we define useful quantities for describing path sets. Then, we  define some generalised travel time models and first arrival times over  path sets. In Remark~\ref{rem:paths} we show how this generalises the standard travel time defined over paths.
For two nodes $x_0,i\in V$, let $\mathbb{P}_{x_0,i}$ be the set of admissible paths $p_{x_0,i}$ 
from $x_0$ to $i$. 
Since the graph $G=(V,E,w)$ is connected, $\mathbb{P}_{x_0,i}$ is non-empty. Let  $P_{x_0,i}\subset \mathbb{P}_{x_0,i}$ denote a non-empty subset of paths from $x_0$ to $i$ and we refer to $P_{x_0,i}$ as a \emph{path set}. We define the \emph{penultimate truncation of a path $p_{x_0,i}\in \mathbb{P}_{x_0,i}$} as a path $p_{x_0,j}$, where $j \in N(i)$ and $p_{x_0,i} = (p_{x_0,j},(j,i))$. Similarly, for a path set $P_{x_0,i}$, we define the \emph{penultimate truncations of $P_{x_0,i}$} as the set $\{ p_{x_0,j} \colon j\in K(P_{x_0,i}) \}$ where $K(P_{x_0,i})\subset N(i)$ such that for every $j\in K(P_{x_0,i})$ there exist a path $p_{x_0,j}$ and a path $p_{x_0,i}\in P_{x_0,i}$ such that $p_{x_0,i} = (p_{x_0,j},(j,i))$. Note, unlike the set $N(i)$ which depends  only on the graph structure, $K(P_{x_0,i})$ depends on the choice of the path set $P_{x_0,i}$. An illustration of a path set and its penultimate truncation is shown in Figure~\ref{fig:pathsets}.

We assume that there exists a formula for a generalised travel time $T(Q)$ for any path set  $Q\subset \mathbb{P}_{x_0,i}$. Some specific examples are introduced below. We define $u_i$ for $i\in V$, as the first arrival travel times over path sets by
\begin{equation}\label{eq:uminpath}
u_i = \min_{x_0\in \partial V}\min_{P_{x_0,i}\subset \mathbb{P}_{x_0,i}} T(P_{x_0,i}).
\end{equation}
For boundary nodes $x_0 \in \partial V$, we set $u_{x_0} = 0$.  
The inner minimisation in \eqref{eq:uminpath} is not over paths $p_{x_0,i}\in\mathbb{P}_{x_0,i}$, but over {\em path sets} $P_{x_0,i}\subset \mathbb{P}_{x_0,i}$.

We define a travel time $T$ over a path set $P_{x_0,i}$ with a local formula over the penultimate truncations of $P_{x_0,i}$. In particular $T(P_{x_0,i})$ is calculated as a function of $T(P^i_{x_0,j})$ with $j \in K(P_{x_0,i})$ where 
$P^i_{x_0,j}=\{p_{x_0,j}\in \mathbb P_{x_0,j} \colon (p_{x_0,j},(j,i))\subset P_{x_0,i}\}$.
By definition $P^i_{x_0,j}$ is also a path set. Since all nodes of a path are distinct by definition, for all $p_{x_0,j}\in P^i_{x_0,j}$ we have $i\notin p_{x_0,j}$.

The models we propose for the travel time $T$ share similarities with the front propagation models 1(i), 1(ii), 1(iii) in Section \ref{sec:frontpropagationmodels} and are specified further below.

\subsubsection{Model 2(i)} 
 Similar to Model 1(i) in \eqref{eq:frontmodel_i}, we define
\begin{equation}\label{eq:arrivalmodel_i}
T(P_{x_0,i}) = \min_{j\in K(P_{x_0,i})} \Big\{T(P^i_{x_0,j}) + \frac{s_i}{w_{j,i}} \Big\}.
\end{equation}
We will  see later that this model is equivalent to the $\ell^\infty$ graph-eikonal equation. 

\subsubsection{Model 2(ii)} 
Similar to Model 1(ii) in \eqref{eq:frontmodel_ii}, we consider
\begin{equation}
T(P_{x_0,i})  = \mu_{x_0,i}
+ \sqrt{ \frac{s_i^2}{z_{x_0,i}} -\sigma_{x_0,i}^2 }\ ,\label{eq:arrivalmodel_ii}
\end{equation}
where 
\[
z_{x_0,i}=\sum_{j\in K(P_{x_0,i})}w^2_{j,i},\qquad 
\mu_{x_0,i}=\frac{1}{z_{x_0,i}} \sum_{j\in K(P_{x_0,i})} w_{j,i}^2T(P^i_{x_0,j})
\]
and
\[\sigma^2_{x_0,i}= \sum_{j\in K(P_{x_0,i})}\left( \frac{ w_{ j,i}^2}{z_{x_0,i}} (T(P^i_{x_0,j}))^2\right)-\mu_{x_0,i}^2.\]
We will  see later that this model is equivalent to the $\ell^2$ graph-eikonal equation.

\subsubsection{Model 2(iii)} 
Similar to Model 1(iii) in \eqref{eq:frontmodel_iii}, we define  
\begin{equation}\label{eq:arrivalmodel_iii}
T(P_{x_0,i})  =\frac{1}{y_{x_0,i}} \sum_{j\in K(P_{x_0,i})} w_{j,i}T(P^i_{x_0,j}) + \frac{s_i}{y_{x_0,i}}
\end{equation}
where $y_{x_0,i}:=\sum_{j\in K(P_{x_0,i})}w_{j,i}$. We will  see later that this model is equivalent to the $\ell^1$ graph-eikonal equation. 

\begin{remark}
Due to the assumption that the graph $G$ is connected and the weights $w_{j,i}$ are positive, there exists a solution to \eqref{eq:uminpath} for all the above choices of the travel time $T$. Clearly first arrival time solutions are well-defined and unique.  However, the minimising path sets are not unique in general.
\end{remark}

\begin{remark}\label{rem:paths}
Consider a singleton path set $P_{x_0,i} = \{p_{x_0,i}\} = \{(x_0 = i_1, \dots, i = i_M)\}$. We observe that the value of $T(P_{x_0,i})$ calculated using models 2(i), 2(ii), or 2(iii), is equal to the following:
\begin{align}\label{eq:traveltime}
\begin{split}
T(\{p_{x_0,i}\})
&= T(\{p_{x_0,i_{M-1}}\}) + \frac{s_{i_M}}{w_{i_{M-1},i_M}}
=T(\{p_{x_0,i_{M-1}}\}) + T(\{(i_{M-1},i_M)\})\\
&=\sum_{m=2}^{M} T(\{(i_{m-1},i_m)\}),
\end{split}
\end{align}
where we used  that the models 2(i), 2(ii), and 2(iii) satisfy 
\begin{equation}\label{eq:edgetraveltime}
T(\{(i_{m-1},i_m)\})=\frac{s_{i_m}}{w_{i_{m-1},i_m}}.
\end{equation}
If we suppose that $w_{i_{m-1},i_m}$  characterises the connectivity between nodes $i_{m-1}$ and $i_m$, and thus $\tfrac{1}{w_{i_{m-1},i_m}}$ is proportional to the travel time, the form of the travel time \eqref{eq:traveltime} can be regarded as a discretisation of $\int_0^1 s(\xi(r)) \|\xi'(r)\|_2\di r$ in \eqref{eq:eikonalmin}. 

Classically, there is a known relationship between the discretisation of problem \eqref{eq:eikonalmin} and the minimisation problem
\begin{equation}\label{eq:M1time}
u_i=\min_{x_0\in \Vbd}\min_{p_{x_0,i}\in \mathbb{P}_{x_0,i}}T(\{p_{x_0,i}\}),
\end{equation}
where $u_{x_0} = 0 $ on boundary nodes $x_0\in \Vbd$. Under the assumption that only singleton sets $P_{x_0,i}=\{p_{x_0,i}\} $ may be considered in \eqref{eq:uminpath}, then \eqref{eq:uminpath} reduces to \eqref{eq:M1time}.

\end{remark}

To understand the behaviour of model 2(i) in \eqref{eq:arrivalmodel_i}, substituting its definition in \eqref{eq:uminpath}, we obtain \eqref{eq:M1time}. Indeed, 
\begin{align*}
u_i &= \min_{x_0\in \partial V}\min_{P_{x_0,i}\subset \mathbb{P}_{x_0,i}} \min_{j\in K(P_{x_0,i})} \Big\{T(P^i_{x_0,j}) + \frac{s_i}{w_{j,i}} \Big\}
\\&= \min_{x_0\in \partial V} \min_{j\in K(\mathbb P_{x_0,i})} \Big\{T(P^i_{x_0,j}) + \frac{s_i}{w_{j,i}} \Big\}
=\min_{x_0\in \Vbd}\min_{p_{x_0,i}\in \mathbb{P}_{x_0,i}}T(\{p_{x_0,i}\}).
\end{align*}
Thus, when using model 2(i), a minimisation over path sets is thus reduced to a minimisation over paths.

\begin{figure}[ht]
    \centering
    \includegraphics[width=\textwidth]{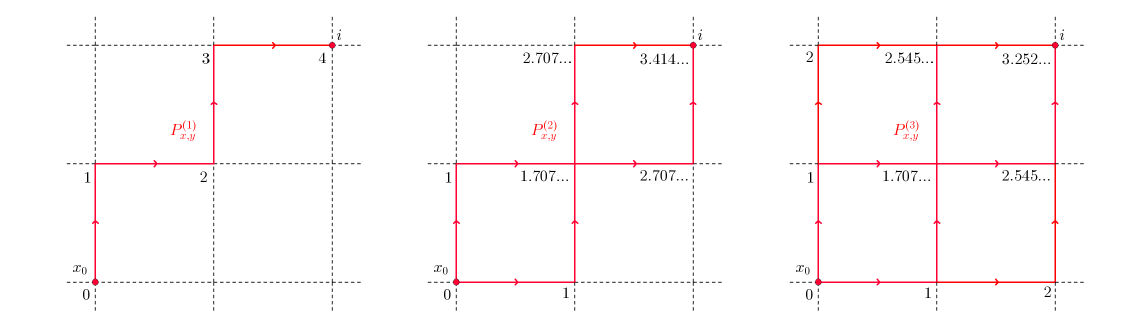}
    \caption{Three different path sets shown in red  on a square grid with $w_{j,i}=1$ and $s_i = 1$ for all nodes. The numbers correspond to the values of the generalised travel time $T(P^{(i)}_{x_0,i})$ for model 2(ii) for each path set.  }
     \label{fig:model3iipaths}
\end{figure}

\begin{figure}[ht]
    \centering
    \includegraphics[width=\textwidth]{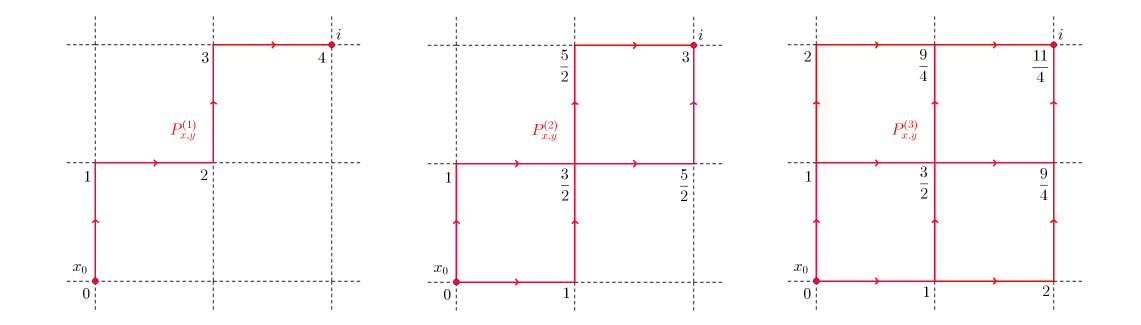}
    \caption{Three different path sets shown in red  on a rectangular grid with $w_{j,i}=1$ and $s_i = 1$ for all nodes. The numbers correspond to the values of the generalised travel time $T(P^{(i)}_{x_0,i})$ for model 2(iii) for each path set.  }
    \label{fig:model3iiipaths}
\end{figure}

To understand the behaviour of models 2(ii) and 2(iii), we calculate the generalised travel time of some simple path sets over the  square grid in two space dimensions with constant unit weights and slowness function; see Figures \ref{fig:model3iipaths} and \ref{fig:model3iiipaths}, respectively. In each case, we calculate the travel times for the three path sets $P^{(1)}_{x_0,i}$, $P^{(2)}_{x_0,i}$ and $P^{(3)}_{x_0,i}$, where $x_0=(0,0)$ and $i=(2,2)$. Let $U$ and $R$ be the paths travelling `up' and `right' from a node to a neighbour on the square grid. We set  $P^{(1)}_{x_0,i}=\{(U,R,U,R)\}$, $P^{(2)}_{x_0,i} = P^{(1)}_{x_0,i} \cup \{(R,U,R,U)\}$ and $P^{(3)}_{x_0,i} = P^{(2)}_{x_0,i} \cup \{(U,U,R,R)$, $(R,R,U,U) \}$, so these path sets have 1, 2 and 4 elements, respectively. We show the generalised travel time for path sets $P^{(1)}_{x_0,i}$, $P^{(2)}_{x_0,i}$ and $P^{(3)}_{x_0,i}$ for models 2(ii) and 2(iii) in  Figures \ref{fig:model3iipaths} and \ref{fig:model3iiipaths}, respectively. Here, the numbers at nodes along the different paths denote the generalised travel time from the origin $x_0$ to the respective nodes. We see that $P^{(3)}_{x_0,i}$ is optimal for model 2(ii) and 2(iii) among $\{P^{(1)}_{x_0,i},P^{(2)}_{x_0,i},P^{(3)}_{x_0,i}\}$ as shown in  Figures \ref{fig:model3iipaths} and \ref{fig:model3iiipaths}. In fact, $P^{(3)}_{x_0,i}$ is an optimal path set for model 2(ii) and 2(iii) among all subsets of $\mathbb{P}_{x_0,i}$ on the square grid.

The properties of minimising path sets are left to future investigation. Heuristically we see that the travel times given by model 2(ii) or 2(iii) are small for path sets that contain short paths or paths which have many cross-overs among themselves (i.e.\ multiple distinct paths pass through common nodes). Such behaviour is observed in Figures~\ref{fig:model3iipaths} and \ref{fig:model3iiipaths}, where the support of the minimizing paths is the rectangular lattice between nodes $x_0$ and $i$.

\begin{remark}
The notion of a minimising path in \eqref{eq:uminpath} also includes the case of a single element of $\partial V$ which corresponds to one label, i.e.\ $\Vbd=\{x_0\}$  in which case
 \[
u_i=\min_{P_{x_0,i}\subset\mathbb{P}_{x_0,i}}T(P_{x_0,i}) .
\]
\end{remark}

\begin{remark}
It is possible that  Gromov's theory (e.g., \cite{Gromov2007}) provides a suitable framework with which to view these constructions. In \cite{Gromov2007}, a metric space $(G,d)$ is endowed with an additional length structure over curves between points in the space. A path-metric space is then defined if $d(x,y)$ for $x,y\in G$ is equal to the shortest length of the curve connecting the $x$ and $y$. This theory applies to the path-distance metric $d$ on a connected graph $G$, defining a length structure by \eqref{eq:arrivalmodel_i}. It will be an interesting future direction of research to see if \eqref{eq:arrivalmodel_ii} and \eqref{eq:arrivalmodel_iii} define length structures and along with a suitable metric $d$ form a path-metric space.
\end{remark}

\subsection{Discrete generalised  eikonal models}
 For $i\in V$, we define one sided edge  derivatives $ \nabla_w^+ u_i \in \mathbb R^{|N(i)|}$ by
\[
    \nabla_w^+ u_i = (w_{j,i}(u_i-u_j)^+)_{j\in N(i)}.
\]
Set,
\begin{equation}\label{eq:pnorm}
\| \nabla_w^+ u_i\|_p=\bl \sum_{j\in N(i)} ({w_{j,i} }(u_i-u_j)^+)^p\br^{1/p} ~~\mbox{for}~~ 1\leq p<\infty,
\end{equation}
 and 
\begin{equation}\label{eq:infinitynorm}
\| \nabla_w^+ u_i\|_\infty= \max_{j\in N(i)} \{{w_{j,i} }(u_i-u_j)^+\}.
\end{equation}

\subsubsection{Model 3(p)}
Motivated by monotone discretisations of the continuum eikonal equation, we consider for any $1\leq p \leq \infty$,
\begin{align}\label{eq:eikonalphelp}
\begin{split}
\| \nabla_w^+ u_i\|_p&=s_i, \quad i\in \mathring V,\\
u_{i}&=0, \quad i\in  \Vbd.
\end{split} 	
\end{align}
 Note that \eqref{eq:eikonalphelp} with $p=2$ is of the same form as the continuum eikonal equation \eqref{eq:eikonal}.
We can rewrite \eqref{eq:eikonalphelp} as 
\begin{align}\label{eq:eikonalp}
\begin{split}
\sum_{j\in N(i)} ({w_{j,i} }(u_i-u_j)^+)^p&=s_i^p, \quad i\in \mathring V,\\
u_{i}&=0, \quad i\in  \Vbd,
\end{split} 	
\end{align}
for $1\leq p<\infty$, and 
\begin{align}\label{eq:eikonalpinfty}
\begin{split}
\max_{j\in N(i)} \{{w_{j,i} }(u_i-u_j)^+\}&=s_i, \quad i\in \mathring V,\\
u_{i}&=0, \quad i\in  \Vbd,
\end{split} 	
\end{align}
for $p=\infty$. 
The models  satisfy a monotonicity condition characteristic  of discrete Hamilton-Jacobi equations (c.f. \cite{DecEllSty11}). Using a monotonicity condition and comparison principles, it has been shown that the boundary value problems 
admit a unique solution and are well-posed, see \cite{DesElmLez13} and \cite[Th. 12]{CalEtt22}. The authors in \cite[Th. 12]{CalEtt22} also construct sub- and supersolutions of the unique solution, resulting in explicit lower and upper bounds of the solution which are both linked to the graph distance.

Note that the $\ell^\infty$ eikonal equation is related to shortest path graph distances which approximate geodesic distances. 
However, this is not the case for the $\ell^p$ eikonal equation with $p$ finite as  interaction between neighbouring nodes are of importance here.

\section{Relations between models}\label{sec:updateformula}
In this section, we investigate relations between the different modelling approaches, that is front propagation, first arrival time and discrete eikonal models, which are introduced in Section \ref{sec:models}. The relationships we prove between the models are summarised in Table~\ref{table:equivalencesmodel} and the proofs are provided in the following sections.

\begin{table}
\centering
 \begin{tabular}{ccccc}\hline\hline
  Front propagation &  &   First arrival & &Discrete eikonal    \\ \hline 
Model 1(i) &   $\Leftrightarrow$ & Model 2(i) & $\Leftrightarrow$ &Model 3($p=\infty$) \\ 
Model 1(ii) & $\Leftrightarrow$ &Model 2(ii) &$\Leftrightarrow$ & Model 3($p=2$) \\ 
Model 1(iii) &    $\Leftrightarrow$ &Model 2(iii) & $\Leftrightarrow$ & Model 3($p=1$) \\ \hline \hline
 \end{tabular}
 \caption{We summarise proved equivalences between the front propagation, arrival time (path and path set) and discrete eikonal models.} \label{table:equivalencesmodel}
 \end{table}

\subsection{Equivalence of front propagation and discrete eikonal models.}
In this section, we show the equivalence of front propagation models \eqref{eq:frontmodel_i}, \eqref{eq:frontmodel_ii}, \eqref{eq:frontmodel_iii} (i.e.\ models 1(i),(ii),(iii)) and discrete eikonal models \eqref{eq:eikonalp} for $p=1, p=2$ and \eqref{eq:eikonalpinfty} for $p=\infty$ (i.e.\ models 3($p=1$), 3($p=2$), 3($p=\infty$)).

\subsubsection{Equivalence of Models 1(i) and 3($p=\infty$)}\label{sec:front(i)eikonal(inf)}

Let $i\in \mathring V$ be given. Hence, there exists $k\in \{1,\ldots,L\}$ such that $i\in V_k$. For this $k$, the definition of sets $V_k, K_{k-1}$ and $C_{k-1}$, and of model 1(i) \eqref{eq:frontmodel_i}, give the value of $u_i$ as 
\begin{align*}
   u_i = U_k = \min_{j\in N(i)\cap K_{k-1}} \left\{u_j + \frac{s_{i}}{w_{j,i}}\right\}, 
\end{align*}
 that is
\begin{align*}
   \max_{j\in N(i)\cap K_{k-1}} \left\{ \frac{w_{j,i}(u_i-u_j) -s_{i}}{w_{j,i}}\right\} =0.
\end{align*}
Since $w_{j,i}>0$ for all edges $(j,i)\in E$, the model is equivalent to
\begin{align*}
   \max_{j\in N(i)\cap K_{k-1}} \left\{ w_{j,i}(u_i-u_j)\right\} -s_{i} =0.
\end{align*}
From minimality of $u_i\in C_{k-1}$, we have $u_j\geq u_i$ for all $j\in V\backslash K_{k-1}$. Recall $w_{j,i}>0$ and $s_i>0$, then extending the set over which the maximum is taken from $N(i)\cap K_{k-1}$ to all of $N(i)$ does not affect the maximum value. Similarly, $K_{k-1}$ necessarily contains at least one point $j$ with $u_j < u_i$, 
therefore replacing $(u_i-u_j)$ with $(u_i - u_j)^+$ does not affect the maximum. This leaves
\begin{align*}
 \max_{j\in N(i)} \left\{ w_{j,i}(u_i-u_j)^+\right\} =s_{i},
\end{align*}
which is precisely \eqref{eq:eikonalpinfty}, i.e.\ model 3($p=\infty$). The counter direction runs exactly the same, with the exception that one must show that $u_i$ is minimal over $C_{k-1}$, however this follows by monotonicity of the construction, as any $j$ with $u_j<u_i$ must belong to $K_{k-1}$ and cannot be in $C_{k-1}$.



\subsubsection{Equivalence of Models 1(ii) and 3($p=2$)}
\label{sec:front(ii)eikonal(2)}
Let $i\in \mathring V$, that is, there exists $k\in \{1,\ldots,L\}$ such that $i\in V_k$.
First, we show that model 3($p=2$) in \eqref{eq:eikonalp}  follows from model 1(ii) in \eqref{eq:frontmodel_ii}. 
For this $k$, the definition of sets $V_k, K_{k-1}$ and $C_{k-1}$, the definition $z_i= \sum_{j\in N(i)\cap K_{k-1}} w_{j,i}^2$, and    \eqref{eq:frontmodel_ii} implies that $u_i$ satisfies
\begin{align*}
\sum_{j\in{N}(i)\cap K_{k-1}} w_{j,i}^2u_i& =  \sum_{j\in{N}(i)\cap K_{k-1}} w_{j,i}^2u_j \\ &\quad + \sqrt{ \bl \sum_{j\in {N}(i)\cap K_{k-1}}  w_{j,i}^2u_j\br^2- z_i \bl\sum_{j\in {N}(i)\cap K_{k-1}} w_{ i,j}^2 u_j^2 -s_i^2\br }. 
\end{align*}
For this, we square both sides of the equality which yields
\begin{align*}
\bl z_iu_i \br^2-2u_i z_i   \sum_{j\in {N}(i)\cap K_{k-1}}  w_{j,i}^2u_j=z_is_i^2- z_i \sum_{j\in {N}(i)\cap K_{k-1}} w_{ i,j}^2 u_j^2. 
\end{align*}
Since $z_i>0$, we obtain
\begin{align}\label{eq:equiv2iihelp}
\sum_{j\in {N}(i)\cap K_{k-1}}  w_{j,i}^2 (u_i-u_j)^2=s_i^2,
\end{align}
From the definition of $K_{k-1}$, the sum can be expanded to the entire neighbourhood $N(i)$, by introducing the maximum with zero,
\begin{align}\label{eq:equiv2iihelp2}
\sum_{j\in {N}(i)}  w_{j,i}^2 ((u_i-u_j)^+)^2=s_i^2,
\end{align}
This is equivalent to model 3($p=2$) in \eqref{eq:eikonalp} 

Next, we start from model 3($p=2$) in \eqref{eq:eikonalp} for $p=2$, or equivalently \eqref{eq:equiv2iihelp}, and show that model 1(ii) in \eqref{eq:frontmodel_ii} follows. Note that \eqref{eq:equiv2iihelp} can be regarded as a quadratic equation in $u_i$ whose solution $u_i$ satisfies
\begin{align*}
u_i=  \frac{1}{ z_i} \bl \sum_{j\in{N}(i)\cap K_{k-1}} w_{j,i}^2u_j  \pm \sqrt{ \bl \sum_{j\in {N}(i)\cap K_{k-1}}  w_{j,i}^2u_j\br^2- z_i \bl\sum_{j\in {N}(i)\cap K_{k-1}} w_{j,i}^2 u_j^2 -s_i^2\br } \br.
\end{align*}
The discriminant is nonnegative due to the existence of a unique real solution to \eqref{eq:eikonalp}. Since 
\begin{align*}
\frac{1}{z_i} \sum_{j\in {N}(i)\cap K_{k-1}} w_{j,i}^2 u_j\leq \max_{j\in {N}(i)\cap K_{k-1}} u_j\leq u_i,
\end{align*}
this implies that the smaller solution contradicts the definition of $i\in V_k$ and the larger solution of the quadratic equation has to be considered, i.e.
\begin{align*}
u_i= \frac{1}{ w_i} \bl \sum_{j\in\tilde{N}(i)} w_{j,i}^2u_j  + \sqrt{ \bl \sum_{j\in \tilde{N}(i)}  w_{j,i}^2u_j\br^2- w_i \bl\sum_{j\in \tilde{N}(i)} w_{j,i}^2 u_j^2 -s_i^2\br } \br,
\end{align*}
which yields \eqref{eq:frontmodel_ii}, that is model 1(ii). $u_i$ is minimal over $C_{k-1}$ by 
 construction, as any $j$ with $u_j<u_i$ must belong to $K_{k-1}$ and thus cannot be in $C_{k-1}$.  

\subsubsection{Equivalence of Models 1(iii) and 3($p=1$)}\label{sec:front(iii)eikonal(1)}
Let $i\in \mathring V$ be given. Hence, there exists $k\in \{1,\ldots,L\}$ such that $i\in V_k$. For this $k$, the definition of sets $V_k, K_{k-1}$ and $C_{k-1}$, and Model 1(iii) in \eqref{eq:frontmodel_iii} show that $u_i$ satisfies
\begin{align*}
u_i  =\frac{1 }{y_i}\bl s_i+\sum_{j\in N(i)\cap K_{k-1}} w_{j,i}u_j \br,
\end{align*}
which is equivalent to model 3($p=1$) in \eqref{eq:eikonalp}  by the definition of $y_i$ and the properties of $i\in V_k$, i.e.
$\sum_{j\in N(i)} w_{j,i}  (u_i-u_j)^+ = s_i$.

\subsubsection{Derivation of Model of type 1 from Model 3($p$) for general $p$}\label{sec:fronteikonal(p)}

We have proved in the previous subsections that there exists a model of type 1 for any model 3($p$), for $p\in\{1,2,\infty\}$. In this subsection we  provide a procedure for deriving such a model of type 1.

For any finite $p\geq 1$, 
the solution $u$ of 3($p$) satisfies \eqref{eq:eikonalp}. Starting from the boundary condition $\partial V$, we initialize the front propagation algorithm. At the $k$th iteration, the following steps are done:
\begin{enumerate}
    \item From $K_{k-1}$ and the graph neighbourhood structure, create $C_{k-1}$.
    \item  By construction of solutions to Model 1, any admissible solution $\tilde u_i$ has to satisfy $\tilde u_i > u_j$ for all $j\in N(i)\cap K_{k-1}$ and $\tilde u_i \leq u_j$ for all $j\in V\backslash K_{k-1}$. 
    To compute the traveltimes $\tilde u_i$ at candidates  $i\in C_{k-1}$, we use \eqref{eq:eikonalp}. Due to the properties of admissible solutions, it is sufficient to restrict the sum in \eqref{eq:eikonalp} to $N(i)\cap K_{k-1}$ instead of $N(i)$. Over this domain, the restriction to the positive part $(\cdot)^+$ may be removed, and the problem is reduced to solving a polynomial equation in $u_i$ (via analytic formulae or numerical solvers).
    As $s_i$ and $w_{j,i}$ are positive, there exists at least one admissible solution $\tilde u_i$. The uniqueness of $\tilde u_i$ follows from contradiction: suppose that \eqref{eq:eikonalp} has two admissible solutions $\bar u_i$ and $\hat u_i$ with $\bar u_i >\hat u_i \geq u_j$ $\forall j\in N(i) \cap K_{k-1}$. 
    Then, 
    $$s_i^p= \sum_{j\in N(i)\cap K_{k-1}} ({w_{j,i} }(\bar u_i-u_j))^p>  \sum_{j\in N(i)\cap K_{k-1}} ({w_{j,i} }(\hat u_i-u_j))^p=s_i^p$$
    This is clearly a contradiction, and thus there is exactly one admissible solution. For all $i\in C_{k-1}$, we denote this admissible solution by $\tilde u_i$ and determine $U_k$ with \eqref{eq:propagationmodel}.
    \item Add all nodes $i\in C_{k-1}$ with  $U_k = \tilde u_i$ into $V_k$, then generate $K_k$.
\end{enumerate}

\subsection{Equivalence of first arrival times over path sets and discrete eikonal models.}
In this section we equate the arrival time model \eqref{eq:uminpath} with travel times \eqref{eq:arrivalmodel_i}, \eqref{eq:arrivalmodel_ii}, \eqref{eq:arrivalmodel_iii} (collectively models 2(i),(ii),(iii)) of Section~\ref{sec:pathset} with the discrete eikonal models, i.e.\ model 3($p=\infty$) in \eqref{eq:eikonalpinfty},  and models 3($p=1$), 3($p=2$) in  \eqref{eq:eikonalp}.

\subsubsection{Equivalence between Models 2(i) and 3$(p=\infty)$ }\label{sec:arrival(i)eikonal(inf)}

Substituting travel time \eqref{eq:arrivalmodel_i} of model 2(i) into \eqref{eq:uminpath} and using the definition of $K(P_{x_0,i})$ for $P_{x_0,i}\subset \mathbb{P}_{x_0,i}$ yields
\begin{align*}
u_i&= \min_{x_0\in \partial V}\min_{P_{x_0,i}\subset \mathbb{P}_{x_0,i}} T(P_{x_0,i})\\
&= \min_{x_0\in \partial V}\min_{P_{x_0,i}\subset \mathbb{P}_{x_0,i}} \min_{j\in K(P_{x_0,i})} \Big(T(P^i_{x_0,j}) + \frac{s_i}{w_{j,i}} \Big)\\
&= \min_{x_0\in \partial V} \min_{K\subset N(i)}  \min_{\{P_{x_0,i}\subset \mathbb{P}_{x_0,i} \colon K(P_{x_0,i}) = K\}} \min_{j\in K} \bl T(P^i_{x_0,j})+\frac{s_i}{w_{j,i}}\br\\
&= \min_{x_0\in \partial V} \min_{K\subset N(i)}\min_{j\in K}  \min_{(P^i_{x_0,j},(j,i))\subset \mathbb{P}_{x_0,i}}  \bl T(P^i_{x_0,j})+\frac{s_i}{w_{j,i}}\br\\
&= \min_{x_0\in \partial V}\min_{j\in N(i)} \min_{(P^i_{x_0,j},(j,i))\subset \mathbb{P}_{x_0,i} }\bl T(P^i_{x_0,j})+\frac{s_i}{w_{j,i}}\br
\end{align*}

Note that $P^i_{x_0,j}$ contains paths between $x_0$ and $j$ not containing node $i$. If we now consider $P_{x_0,j}\subset \mathbb{P}_{x_0,j}$, then there may be a path from $x_0$ to $j$ via $i$ in $P_{x_0,j}$, but it is not a minimiser. To see that a path $p_{x_0,j}$ with $i\in p_{x_0,j}$ is indeed not a minimiser, we consider $p_{x_0,j}=(i_1=x_0,\ldots,i_k=i,\ldots,i_M = j)$ for some $M\in \N$ and $1<k<M$, implying that $i_{k-1}\in N(i)$ and hence  $p_{x_0,i}=(i_1=x_0,\ldots,i_{k-1}, i_k=i)\in (P^i_{x_0,\tilde j},(\tilde j,i))\subset \mathbb{P}_{x_0,i}$ for $\tilde j=i_{k-1}\in N(i)$ and some path set $P^i_{x_0,\tilde j}\subset \mathbb P_{x_0,\tilde j}$. As the travel time is nonnegative on every edge by \eqref{eq:edgetraveltime}, the travel time is monotone over increasing path length and we have $T(p_{x_0,i_{k-1}})<T(p_{x_0,j})$ with $i_{k-1},j\in N(i)$, implying that $p_{x_0,j}$ with $i\in p_{x_0,j}$ cannot be a minimiser. Hence, we write
\begin{align*}
    u_i&= \min_{x_0\in \partial V}\min_{j\in N(i)} \min_{P_{x_0,j}\subset \mathbb{P}_{x_0,j} }\bl T(P_{x_0,j})+\frac{s_i}{w_{j,i}}\br\\
&= \min_{j\in N(i)}\bl\bl\min_{x_0\in \partial V} \min_{P^i_{x_0,j}\subset \mathbb{P}_{x_0,j}}T(P_{x_0,j})\br +\frac{s_i}{w_{j,i}} \br=\min_{j\in N(i)}\bl u_j +\frac{s_i}{w_{j,i}}\br.
\end{align*}
We move $u_i$ to the right hand side, and use that $\min(x) = - \max(-x)$, so that we obtain
\begin{align*}
0=\max_{j\in N(i)} \bl u_i-u_j - \frac{s_{i}}{w_{j,i}}\br=\max_{j\in N(i)} \bl\frac{w_{j,i} \bl u_i-u_j \br-s_{i}}{w_{j,i}}\br.
\end{align*}
Due to the positivity of $w_{ij}$, this is equivalent to $\max_{j\in N(i)} \bl w_{j,i} \bl u_i-u_j \br-s_{i}\br=0$, 
and as $u_i\geq u_j$, this yields
$\max_{j\in N(i)} \bl w_{j,i} \bl u_i-u_j \br^+\br=s_{i}$, 
that is, we obtain model 3($p=\infty$) in \eqref{eq:eikonalpinfty}.

\subsubsection{Equivalence between Models 2(ii) and 3$(p=2)$ }\label{sec:arrival(ii)eikonal(2)}

Starting with \eqref{eq:uminpath} and considering travel time of model 2(ii) in \eqref{eq:arrivalmodel_ii}  yields
\begin{align*}
    u_i&= \min_{x_0\in \partial V}\min_{P_{x_0,i}\subset \mathbb{P}_{x_0,i}} T(P_{x_0,i})\\
    &= \min_{x_0\in \partial V}\min_{P_{x_0,i}\subset \mathbb{P}_{x_0,i}} \left( \frac{1}{z_{x_0,i}} \sum_{j\in K(P_{x_0,i})} w_{j,i}^2T(P^i_{x_0,j})\right.
\\&\quad\left.+\frac{1}{z_{x_0,i}} \sqrt{ \bl  \sum_{j\in K(P_{x_0,i})} w_{j,i}^2T(P^i_{x_0,j}) \br^2+ z_{x_0,i}s_i^2  -z_{x_0,i}\sum_{j\in K(P_{x_0,i})}  w_{ j,i}^2 (T(P^i_{x_0,j}))^2 }\right)
\end{align*}
where
$z_{x_0,i}=\sum_{j\in K(P_{x_0,i})}w^2_{j,i}$.
We can write $u_i$ as
\begin{align*}
u_i&=\min_{x_0\in \partial V}\min_{K\subset N(i)} \min_{\{P_{x_0,i}\subset \mathbb{P}_{x_0,i}\colon K(P_{x_0,i})=K\}} \left(\frac{1}{z_K} \sum_{j\in K} w_{j,i}^2T(P^i_{x_0,j})\right.
\\&\qquad\qquad \left.+\frac{1}{z_K} \sqrt{ \bl  \sum_{j\in K} w_{j,i}^2T(P^i_{x_0,j}) \br^2+ z_K s_i^2  -z_K\sum_{j\in K}  w_{ j,i}^2 (T(P^i_{x_0,j}))^2 }\right),
\end{align*}
where 
$z_K=\sum_{j\in K}w^2_{j,i}$.
Since $T(P^i_{x_0,j})$ is the only term depending on $x_0\in \partial V$ and $P^i_{x_0,j}$ satisfying $P_{x_0,i}=(P^i_{x_0,j},(j,i))\subset \mathbb{P}_{x_0,i}$ with $j\in K(P_{x_0,i})$, we may pull the minimisation with respect to these parameters inside the expression, and replace the minimisation with respect to $P_{x_0,i}=(P^i_{x_0,j},(j,i))\subset \mathbb{P}_{x_0,i}$ with $j\in K(P_{x_0,i})$ by $P_{x_0,j}\subset \mathbb P_{x_0,j}$  as in Section~\ref{sec:arrival(i)eikonal(inf)}. This yields
\begin{align*}
    u_i&=\min_{K\subset N(i)} \left( \frac{1}{z_K} \sum_{j\in K} w_{j,i}^2u_j
+\frac{1}{z_K} \sqrt{ \bl  \sum_{j\in K} w_{j,i}^2u_j \br^2+ z_K s_i^2  -z_K\sum_{j\in K}  w_{ j,i}^2 u_j^2 }\right)
\end{align*}
where $u_j= \min_{x_0\in \partial V}\min_{P_{x_0,j}\subset \mathbb{P}_{x_0,j}} T(P_{x_0,j})$ by definition.
Moving  $u_i$ to the right-hand-side and using $\min(x)=-\max(-x)$ provides 
\begin{align*}
    0=\max_{K\subset N(i)}\left( \frac{1}{z_K}\sum_{j\in K}w^2_{j,i} (u_i - u_j) -\frac{1}{z_K} \sqrt{ \bl  \sum_{j\in K} w_{j,i}^2u_j \br^2+ z_K s_i^2  -z_K\sum_{j\in K}  w_{ j,i}^2 u_j^2 }\right).
\end{align*}
To achieve that the expression vanishes, we require that the first term is nonnegative which is equivalent to $K\subset N(i)$ such that $u_j\leq u_i$ for all $j\in K$. Note that the first term is maximal for the set  $\{j \in N(i) \colon u_j \leq u_i\}$ and the magnitude of the second term decreases as the size of the set $K$ increases. 
Hence, the maximiser $K$ with $K=\{j \in N(i) \colon u_j \leq u_i\}$ satisfies
\begin{align*}
     z_Ku_i-\sum_{j\in K}w^2_{j,i} u_j = \sqrt{ \bl  \sum_{j\in K} w_{j,i}^2u_j \br^2+ z_K s_i^2  -z_K\sum_{j\in K}  w_{ j,i}^2 u_j^2 }.
\end{align*}
Squaring both sides and dividing by $z_K$ yields
\begin{align*}
     z_Ku_i^2-2u_i \sum_{j\in K}w^2_{j,i} u_j   =   s_i^2  -\sum_{j\in K}  w_{ j,i}^2 u_j^2, 
\end{align*}
i.e.
\[
s_i^2 = \sum_{j\in K}  w_{j,i}^2 (u_i-u_j)^2 = \sum_{\substack{j\in N(i) \colon u_j\leq u_i}}  w_{j,i}^2 (u_i-u_j)^2 = \sum_{j\in N(i)}w_{j,i}^2 ((u_i-u_j)^+)^2,
\]
that is model 3($p=2$) in \eqref{eq:eikonalp}.

\subsubsection{Equivalence between Models 2(iii) and 3$(p=1)$ }\label{sec:arrival(iii)eikonal(1)}

We begin by using the first arrival model \eqref{eq:uminpath} with travel time $T$ given as in model 2(iii) by \eqref{eq:arrivalmodel_iii} which yields
\begin{align*}
u_i &= \min_{x_0\in \partial V}\min_{P_{x_0,i}\subset \mathbb{P}_{x_0,i}} T(P_{x_0,i})\\
&= \min_{x_0\in \partial V}\min_{P_{x_0,i}\subset \mathbb{P}_{x_0,i}}\frac{1}{\sum_{j\in K(P_{x_0,i})}w_{j,i}} \bl \sum_{j\in K(P_{x_0,i})} w_{j,i}T(P^i_{x_0,j}) + s_i\br\\
&=\min_{K\subset  N(i)} \min_{x_0\in \partial V} \min_{\{P_{x_0,i}\subset \mathbb{P}_{x_0,i}\colon  K(P_{x_0,i})=K\}}\frac{1}{\sum_{j\in K}w_{j,i}} \bl \sum_{j\in K} w_{j,i}T(P^i_{x_0,j}) + s_i\br\\
&=  \min_{K\subset  N(i)}
\frac{1}{\sum_{j\in K}w_{j,i}} \bl \sum_{j\in K} w_{j,i} \min_{x_0\in \partial V}\min_{P_{x_0,j}\subset \mathbb{P}_{x_0,j}} T(P_{x_0,j}) + s_i\br\\
&= \min_{K\subset N(i)} \frac{1}{\sum_{j\in K}w_{j,i}}\bl \sum_{j\in K} w_{j,i}u_j+ s_i \br,
\end{align*}
where we can use a similar argument as in Section \ref{sec:arrival(i)eikonal(inf)} in the fourth equality to consider the sets $P_{x_0,j}\subset \mathbb P_{x_0,j}$ instead of the sets $P_{x_0,i}=(P^i_{x_0,j},(j,i))\subset \mathbb{P}_{x_0,i}$ with $j\in K(P_{x_0,i})$.
Then we rearrange the equation resulting in 
\begin{align*}
    \min_{K\subset N(i)} \frac{1}{\sum_{j\in K}w_{j,i}}\bl \sum_{j\in K} w_{j,i}(u_j-u_i)+ s_i\br=0,
\end{align*}
and as $\sum_{j\in K}w_{j,i}>0$, we obtain
\begin{align*}
 s_i & = -\min_{K\subset N(i)} \bl  \sum_{j\in K} w_{j,i}  (u_j-u_i) \br =\max_{K\subset N(i)} \sum_{j\in K} w_{j,i}  (u_i -  u_j ).
\end{align*}
If $u_j\leq u_i$ then the summand is positive and therefore the maximiser over $K\subset N(i)$ is the set  $\{j \in N(i) \colon u_j \leq u_i\}$. Hence we arrive at
\[
 s_i = \sum_{\substack{j\in N(i)\colon  u_j \leq u_i}} w_{j,i}  (u_i -  u_j ) =\sum_{j\in N(i)} w_{j,i}  (u_i -  u_j )^+,
\]
that is, model 3($p=1$) in \eqref{eq:eikonalp}.

\section{Applications}

In this section, we use information propagation in two applications: trust networks and semi-supervised learning.

\subsection{Applications to trust networks\label{sec:trustnetworks}}

In this section, we apply information propagation to a trust network, which is a weighted directed graph, with nodes being users of a social network. The edges and their trust weights are reviews of trust among users, for example, $\omega_{i,j}$ implies that $i$ \emph{trusts} $j$ with a rating $\omega_{i,j}$. The neighbourhood structure is therefore built around \emph{trusting nodes} and their neighbours of \emph{trusted nodes}. It is a directed relationship, as trust may not be reciprocated (i.e.\ $(i,j)$ may exist, but $(j,i)$ does not), and may not be symmetric ($\omega_{i,j}\neq \omega_{j,i}$).

An example of a trust network is the advocato dataset \cite{DanMarPao09}, specifically \texttt{soc-advocato} obtained from \cite{AhmRos}. The dataset is a snapshot in time of a social network comprised of around 5000 software developers, with four categorical weightings of trust based on  users' perceived contributions to open software and programming skills. These weightings have been numerically equated, though somewhat arbitrarily  \cite{DanMarPao09}, to four values $\omega_{i,j} \in (0.4,0.6,0.8,1.0)$ where larger values correspond to larger levels of trust. Structural information of the graph is found in \cite{AhmRos}. For this experiment we will investigate only the largest connected component of this graph, which contains 5167 nodes (of which 4017 are trusting, and 4428 are trusted) and 47337 edges.

Consider the application of a software team searching for a new collaborator from the network. The team seeks a notion of trustworthiness for each candidate collaborator. One can assess the level of trust in each candidate by the team by using the information from the trust network. A seemingly simple way to do this is to aggregate the trust given to them directly by other users (i.e. the weights from users to the candidate). Unfortunately, a common subversion of this method is a  Sybil attack \cite{Dou02,LuSheZha14,AlQ_etal17}. In its simplest form, a candidate creates a network of artificial community members, colloquially called ``Sybils", who have high trust with each other and with the candidate. This will boost the candidates aggregate trust. Instead of a neighbourhood-only measurement of trust, we propose using information propagation of distrust to provide a ranking candidates from the perspective of the team in a way that is resistant to Sybil attack.

We perform the propagation of distrust by setting boundary nodes $\partial V$ as the software team. We then set $s=1$ at all nodes. We define the distrust weights $w_{j,i} = \frac{1}{\omega_{j,i}}$, a reciprocal of the trust weights. Candidates are selected on the network, and we calculate the (first arrival) times for the information to propagate to all candidates from the team members by using a propagation model. The candidates with greater arrival time are less trustworthy according to the model. This method of measuring distrust accounts for both the degree of separation between the team and the candidates, as well as the trust of each review along such paths. It is resistant to Sybil attacks, as ``Sybils" form a (largely) disconnected cluster around a candidate, and so have little or no effect on path structures between the team and the candidate. 

\begin{figure}
   \centering
   \includegraphics[height=5.35cm]{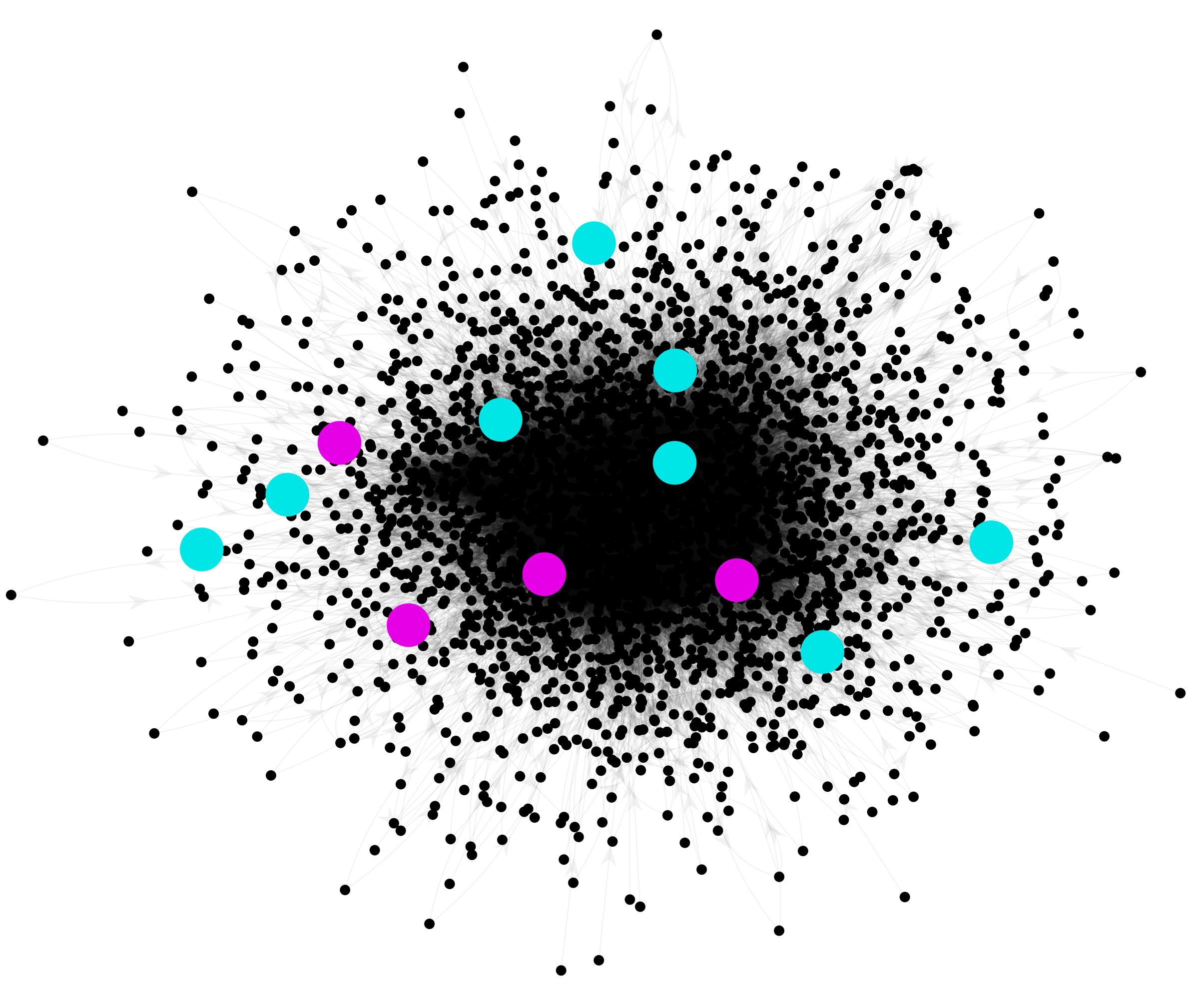}
   \includegraphics[height=5.35cm]{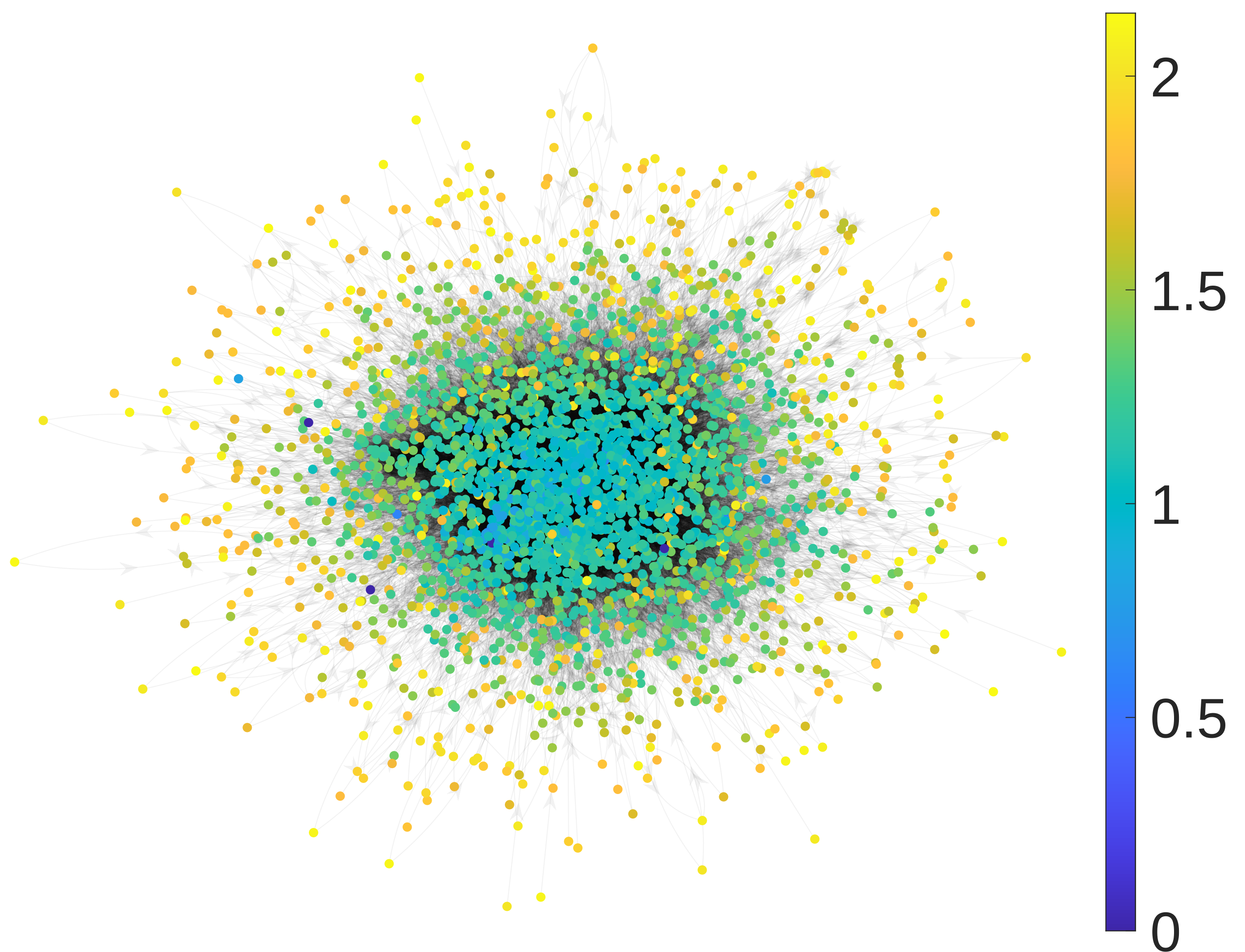}
\caption{Result of the distrust propagation from a four-member software team, to eight candidates. Edge arrows indicates direction of trust. The left panel shows the software team (magenta) and candidates (cyan). The right panel shows the solved travel time field using model $p=1$, with node colour indicating the level of distrust of this community member by the software team.}
\label{fig:l1trustrank}
\end{figure}

The experiment configuration is shown overlaying a relevant portion of the network in the left panel of Figure \ref{fig:l1trustrank}, we randomly select both a four-member team in magenta, and eight candidates in cyan that we label $A$-$H$. For illustration, the travel times (level of distrust) of the displayed nodes according to the propagation model with $p=1$ are given in the right panel of Figure \ref{fig:l1trustrank}. We perform two control experiments. In the first control experiment (Ctrl), we directly use the network of \texttt{soc-advocato}. In a second experiment (GSyb), we modify the network by adding a fully connected Sybil cluster (of size 50) to candidate $G$. The members of this cluster are given by the highest trust weighting 1 from each other and the candidate $G$, and vice versa. For each experiment, we use different information propagation models, and the neighbour-averaging approach to calculate a trust ranking of the eight candidates. The results of the experiments are given in Table~\ref{tab:trustrank}. Candidates $A$-$H$ are alphabetically assigned by the order of the first column. The first three columns show ranks given by propagation model for $p=1,2,\infty$ which is based on globally averaging distrust. The final two columns show ranks for candidate $j$ based on locally averaging distrust $w_{i,j}$ over $i$ such that $j\in N(i)$. 

First, we look at  experiment (Ctrl) in Table \ref{tab:trustrank}. Comparing the propagation-based and neighbour-based ranks, we offer an interpretation of some interesting candidates. As candidate $A$ is deemed trustworthy across all methods, this implies that both locally and globally $A$ is a trustworthy candidate. Candidate $H$, on the other hand, is deemed trustworthy locally but suffers globally which indicates that there is an overall untrustworthy pathway of reviewers between the team and $H$. Candidate $D$ shows a difference in ranking between travel time models $p=1,2$ and $p=\infty$ which implies that although the most trustworthy review path from the team to $D$ is not very trustworthy ($p=\infty$ ranks $D$ eighth), there are many similarly trusted pathways from the team to $D$ ($p=1,2$ rank $D$ fourth and fifth). In this way, the models for $p<\infty$ encode a concept of confidence over the network uncertainty into their travel time. In general, we see that for $p=1$ and $p=2$, candidates enjoy similar rankings, while the $p=\infty$ model tends to group candidates together, as for this model the travel times can only take more restrictive discrete values. These preliminary results suggest that $p < \infty$ approaches provide solutions richer in information from the network, and may be more robust in discrete settings, and so we advocate further investigation of their use in graph-based algorithms where $p=\infty$ may be the state-of-art.

Finally, the key result of comparing between the experiments (Ctrl) and (GSyb) is that the travel time based ranking did not change between the experiments, whereas the neighbour-averaged  distrust of candidate $G$ reduced from $1.333$ in (Ctrl) to $1.029$ in (GSyb), thus increasing their rank from 7 to 4. This provides concrete evidence to the susceptibility of neighbourhood-based approaches, while arrival time approaches are completely resilient to this form of attack. 

\begin{table}
\begin{tabular}{cccccc} \hline \hline
 Trust Rank   & $p=1$&   $p=2$  & $p=\infty$ & Neighbour & Neighbour \\  
 (1 = most trust) & Ctrl \& GSyb  & Ctrl \& GSyb & Ctrl \& GSyb & Ctrl & GSyb \\ \hline
 1 & $A$ (1.170) &  $A$ (1.728) & $B$ (2.0) & $A,D,H$ (1.0) & $A,D,H$ (1.0) \\ 
 2 & $B$ (1.360) &  $B$ (2.001) & $A,C,E$ (2.2) & $A,D,H$ (1.0) & $A,D,H$ (1.0)\\
 3 & $C$ (1.371) &  $C$ (2.019) & $A,C,E$ (2.2) & $A,D,H$ (1.0) & $A,D,H$ (1.0)\\
 4 & $D$ (1.573) &  $E$ (2.084) & $A,C,E$ (2.2) & $B$ (1.167) & $\mathbf{G (1.029)}$\\
 5 & $E$ (1.784) &  $D$ (2.341) & $F,G$ (2.4)   & $C,F$ (1.25) & $B$ (1.167)\\
 6 & $F$ (2.005) &  $F$ (2.354) & $F,G$ (2.4)   & $C,F$ (1.25) & $C,F$ (1.25)\\
 7 & $G$ (2.047) &  $G$ (2.358) & $H$ (2.6)     & $\mathbf{G (1.333)}$  & $C,F$ (1.25) \\
 8 & $H$ (2.148) &  $H$ (2.574) & $D$ (2.8)     & $E$ (1.458)  & $E$ (1.458)\\ \hline \hline
\end{tabular}
\caption{Ranking of trust of candidates $A$-$H$, for two experiments: a control experiment (Ctrl), and an experiment with a cluster of 50 Sybils around candidate $G$ (GSyb). Candidates $A$-$H$ are alphabetically assigned by the order of the first column. The columns give trust rankings from different information propagation models ($p=1,2,\infty$), or from using the average of neighbourhood distrust (Neighbour). The measure of absolute distrust of the candidate is given in brackets: for the first three columns this is the travel-time, in the final two columns this is the averaged distrust over the neighbourhood of the candidate. }
\label{tab:trustrank}
\end{table}

\subsection{Applications to label propagation/semi-supervised learning}\label{SemiSuper}

In this section, we consider an application to a semi-supervised learning approach to label propagation.
The model consists of attaching $L>1$ labels to $n>1$ sets of features
$f_j\in \mathcal F_j,~j=1,2,\ldots,n$, where $\mathcal F_j=\{\mathcal F^i_j\in \mathbb F_i\}_{j=1}^m~i=1,2,\ldots,m$ and $\mathbb F_i~~\mbox{is either}~\mathbb R~~\mbox{or}~~\mathbb B,~\mathbb B =\{0,1\}$. The first step consists of assigning weights $w_{i,j}\ge0$ whose reciprocal measures the distance between features $f_i$ and $f_j$. If the distance between features is sufficiently large according to some criterion then the weight is set to zero.  From this we obtain the graph with vertices $V=\{1,2,\ldots, n\}$ and edges $E\in V^2$ consisting of the  pairs satisfying $w_{i,j}>0.$ 
We assume there is a set of nodes $\partial V_\ell$ for each category $\ell=1,\dots, L$ so that label$(i) = \ell$ for all $i\in\partial V_\ell$, that is, a set where the classification is known. Our semi-supervised learning task is to provide all unlabelled nodes in $V\setminus \{\cup_\ell \partial V_\ell\}$ with a label. 
The front propagation semi-supervised learning model is to assign  
\begin{equation} \label{eq:lbl_strategy}
\text{label}(i) = \{\ell \ | \ u^{(\ell)}_i = \underset{k=1,\dots,L}{\min} u^{(k)}_i\}
\end{equation}
for any $i \in V\setminus \{\cup_{\ell} \partial V_\ell\}$,
where $u^{(k)}_i$ is the solution to a discrete eikonal equation \eqref{eq:eikonalphelp} on a weighted graph for some $p\in [1,+\infty]$, with boundary $u_i=0$ for $i\in\partial V_k$. We assume for this model that the slowness function $s\equiv 1$. In this way, $i$ is assigned the label $\ell$ if the smallest travel time is found between $\partial V_\ell$ and $i$ among all sets of labels. This model requires to solve  the discrete eikonal problem per label category, which can be performed independently in parallel to each other. 
For each of the following experiments, we carry out  20 simulations  with differing random choices of known initial labels. We present the average (and standard deviation) of the classification accuracy over these 20 simulations.  The labelling accuracy is calculated as the percentage of unlabelled nodes that are correctly classified.

\subsubsection{High-dimensional two moons problem}
\begin{figure}[htb]
\centering 
\includegraphics[width=\textwidth]{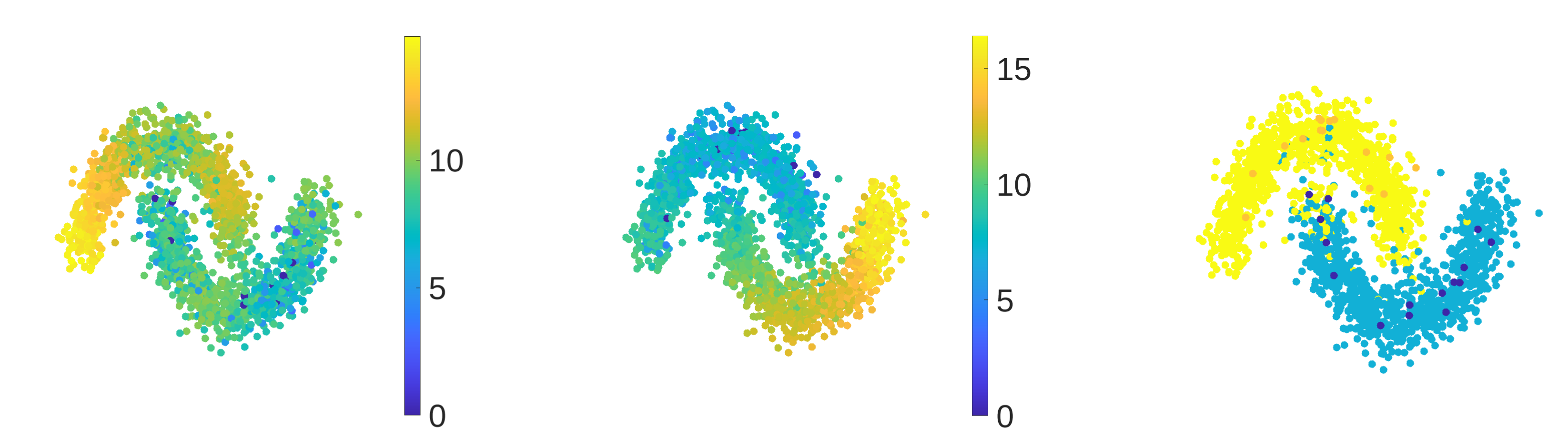}
\caption{Example travel time fields and classification for two moons problem, projected into two dimensions. The left and centre panels show the travel time field for label 1 and 2 respectively. The right panel shows the resulting classification with predicted label 1 (blue) and predicted label 2 (yellow) solved with initially known labels 1 (orange), and 2 (dark blue). In this example, the accuracy was 94.7\%.}
\label{fig:twomoons}
\end{figure}

We follow the construction of the two moons problem for classification as in \cite{BuhHei09,BerFle16}. The feature vectors here are taken to be the spatial coordinates of $n$ nodes in $\mathbb{R}^{m},$~ i.e.~$\mathbb F_i= \mathbb R, ~\forall i.$  The construction is formed by considering two planar half circles of radius 1. One is centred at the origin and  the other is rotated by $\pi$ and centred $(1,0.5)$.  We take $n=2000$ points on these  initial  planar circles  and then embed them in $\mathbb R^{100}$ by adding 
uniform Gaussian noise $N(0,0.02I_{100})$ where $I_{100}$ is the identity matrix  in $\mathbb R^{100}$.
We define a classification problem by giving points on each initial circle  a different binary  label; for visualisation we project back onto the plane as seen in Figure~\ref{fig:twomoons}. 
We proceed again as in \cite{BuhHei09,BerFle16} by calculating distances between pairs of points in $\mathbb R^{100}$ and then setting all weights $w_{i,j}=0$ unless point $j$ is within the 10 nearest neighbours of point $i$. The non-zero weights are then set according to the weight function of \cite{PerZel04}; a squared exponential function of distance, weighted by a local scaling $d_{10}(x_i)=\|x_i-x_{j(i,10)}\|$, where $j(i,10)$ is the 10$^{\text{th}}$ nearest neighbour of $i$ (see Table~\ref{tab:twomoon_accuracy}). 
We perform each of the experiments by choosing at random  15 nodes per moon  to have  known labels. The illustration of the travel time-based classification is given in Figure~\ref{fig:twomoons}. The accuracy results are given in Table~\ref{tab:twomoon_accuracy}. Here we observe high accuracy, with all choices of eikonal model comparable to experiments of unsupervised clustering in \cite{BerFle16} with near optimal parameter choices. Our method has no tuneable parameters, though the experiment suggests best performance for $p=1$. 

\begin{table}
    \centering
    \begin{tabular}{ccc}\hline \hline
       $w_{i,j}$  & Eikonal model & Two moons accuracy \% \\ \hline
                            &$p=1$ & 92.7 (3.81) \\ 
$\exp(-\frac{\|x_i - x_j\|^2}{\sqrt{d_{10}(x_i)d_{10}(x_j)}})$  &$p=2$      & 92.0 (2.80)\\ 
                            &$p=\infty$ & 89.5 (2.96) \\ \hline \hline
    \end{tabular}
    \caption{Mean (standard deviation) of classification for the two moons example.}
    \label{tab:twomoon_accuracy}
\end{table}
\subsubsection{Text classification dataset}
We demonstrate the performance on the standard Cora and CiteSeer document classification datasets \cite{Sen_etal08}. In both cases, the graph nodes correspond to  journal articles, and links between them are obtained from  citations, forming a directed graph. The featue vectors are binary valued of length  1433 (Cora, ~i.e. $\mathbb F_i=\mathbb B,~\forall i$) and  3703 (CiteSeer~ $\mathbb F_i=\mathbb B,~\forall i$), based on whether or not  the article contained specific words from a unique dictionary. Following  \cite{CohSalYan16,pre:kipfWel16}, we symmetrise the adjacency matrix for each citation link.  We benchmark with the resulting largest connected component of each dataset. The resulting graphs have  2485 nodes and 5069 edges (Cora), and 2110 nodes and 3694 edges (CiteSeer). The reference did not provide suggestion for the graph weights, thus some naive choices were taken, based on the $\ell^2$-norm over binary vectors (see Table~\ref{tab:text_data_accuracy}).  
There are  $L=7$ (Cora) and  $L=6$ (CiteSeer) labels respectively for each dataset, representing journal categories that we wish to classify. We take 20 labels from each category. The classification accuracy experiments for the different data sets were applied to 20 random seeds, and we display the results of the eikonal models for $p\in\{1,2,\infty\}$. We assume for this application that the slowness function $s\equiv 1$. The results are shown in Table~\ref{tab:text_data_accuracy}. Performance is robust across seeds and eikonal models chosen. The experiment suggests best performance at $p=1$. The exponential based weighted graphs outperform the reciprocal distance based weights, and have less variation due to random seeding. For this graph, $d_{\max}$ was relatively constant and did not aid performance. We did not optimise the constants appearing in the weight functions, and the algorithms performed similarly across an order of magnitude.
Several approaches have applied to these data sets in \cite{CohSalYan16}. Here comparisons are qualitative, as different methods (e.g., \cite{Joa99,BelNiySin06,ColMobRatWes12}) use differing levels of information. On these data sets, the front propagation approach performs comparably to Planetoid-T and Planetoid-I, the flagship methods of \cite{CohSalYan16}.

\begin{table}
    \centering
    \begin{tabular}{cccc} \hline \hline
       $w_{i,j}$  & Eikonal model &   Cora accuracy & CiteSeer accuracy \\ \hline
     
                            &$p=1$       & 69.0 (7.49) & 64.3 (1.64) \\
$ 1/\|x_i - x_j\|_{\ell^2}$ &$p=2$      & 68.9 (6.86) & 62.6 (1.87) \\
                            &$p=\infty$  & 68.1 (3.86)& 61.0 (2.26)\\ \hline
    
        
                                     &$p=1$      & 72.4  (1.58) & 64.3 (1.91) \\
$\exp(-\frac{\|x_i - x_j\|^2_{\ell^2}}{500})$&$p=2$       & 71.8 (1.88) & 62.5 (2.12) \\
                                     &$p=\infty$  &69.2 (2.50)& 60.8 (2.25)\\ \hline
                                      &$p=1$      & 72.4  (1.56) &  64.3 (2.06) \\
$\exp(-\frac{\|x_i - x_j\|^2_{\ell^2}}{100\sqrt{d_{\max}(x_i)d_{\max}(x_j)}})$&$p=2$       & 71.7 (1.91) & 62.5 (2.08) \\
                                     &$p=\infty$ & 69.0 (2.42)& 60.7 (2.22)\\ \hline \hline
    \end{tabular}
    \caption{Mean (standard deviation) of classification accuracy given as percentages, for the examples using different choices of weights. The function $d_{\max}(x)$ is the Euclidean distance from $x_i$ to its furthest neighbour.}
    \label{tab:text_data_accuracy}
\end{table}

\section{Conclusion}\label{Conclusion}
In this paper we proposed some models for information propagation on graphs. Underlying components of the models include a subset of nodes forming  the  initial source of information, the arrival times of  information  and the  laws governing the transfer of information to nodes from their neighbours. The models are collected into three viewpoints: an information wavefront hitting time, an optimal travel time over sets of paths, and a local equation for the time to receive information at a node given the times to receive information at its neighbours.  
We showed equivalences between these different views, as summarised in Table~1. 
In this framework we provide examples such as  a generalisation of classical equivalence between optimal paths and Djikstra's algorithm \cite{Dji59}. 
We applied our models to a social network dataset \texttt{soc-advogato} \cite{AhmRos}, where directed edges are weighted by trust. Propagation of wavefronts from a group of nodes over such weighted networks define a notion of (dis)trust of this group on all other nodes defined by the travel times to other nodes. This notion of trust is robust to local Sybil attacks \cite{Dou02}. More generally, our models could be used as a back-end to replace path-length or distance calculations in other cybersecurity strategies \cite{AlQ_etal17}, as qualitatively the $p<\infty$ approaches displayed better solution properties than $p=\infty$. 
Extending the work of  \cite{Tou14}, we applied these models to label propagation in a semi-supervised learning application. The eikonal-based classification algorithm obtains comparable performance to clustering algorithms with two labels (e.g., \cite{BerFle16}), and with simple choices of weight functions, it achieves comparable performance to machine-learning methods that learn graph embeddings (e.g., \cite{CohSalYan16}) without any tuning or training. 
While graph Laplacian methods are often used to model information propagation on networks (e.g.,  \cite{elmdestou17,Elm15,FaxMurOlf07}), the eikonal approach can also be applied and
encapsulates control problems (using $s$ or $w$ as controllers). Procedures based on front propagation offer adjoint equations at no additional cost, leading to very efficient methods for inverse problems in these settings \cite{DecEllSty11,DunEll19}.

\section*{Acknowledgements}
ORAD would like to acknowledge the support of Schmidt Sciences, LLC, the  National Science Foundation (Grant No. AGS-1835860), the Cisco Foundation, and the Office of Naval Research (Grant No. N00014-23-1-2654). LMK acknowledges support from the the Warwick Research Development Fund through the project `Using Partial Differential Equations Techniques to Analyse Data-Rich Phenomena', the European Union Horizon 2020 research and innovation programmes under the Marie Sk\l odowska-Curie grant 
agreement No.\ 777826 (NoMADS), the Cantab Capital Institute for the Mathematics of Information and Magdalene College, Cambridge 
(Nevile Research Fellowship).

\bibliographystyle{siam}      
\bibliography{references}

\end{document}